\documentclass[12pt,reqno]{amsart}
\ifx\pdfsuppressptexinfo\undefined\else\pdfsuppressptexinfo=-1\fi
\usepackage{amsmath,amsthm,amssymb}
\usepackage{enumitem}
\usepackage{needspace}
\usepackage[hmargin=1.15in,vmargin=1.15in]{geometry}
\usepackage[colorlinks,linkcolor=blue,citecolor=blue,urlcolor=blue]{hyperref}
\usepackage{bookmark}
\setlist[enumerate,1]{label=(\arabic*),ref=\arabic*,font=\normalfont}
\numberwithin{equation}{section}
\newcommand{\F}{\mathbb F}
\newcommand{\HH}{\mathcal H}
\newcommand{\LL}{\mathcal L}

\DeclareMathOperator{\Aut}{Aut}
\DeclareMathOperator{\Gal}{Gal}
\DeclareMathOperator{\Fix}{Fix}
\DeclareMathOperator{\Fr}{Fr}
\DeclareMathOperator{\Tr}{Tr}
\DeclareMathOperator{\ord}{ord}
\DeclareMathOperator{\Nrm}{N}
\DeclareMathOperator{\Cen}{C}
\DeclareMathOperator{\SL}{SL}
\DeclareMathOperator{\PGL}{PGL}
\providecommand{\X}{\mathcal X}

\theoremstyle{plain}
\newtheorem{theorem}{Theorem}[section]
\newtheorem{proposition}[theorem]{Proposition}
\newtheorem{lemma}[theorem]{Lemma}
\newtheorem{corollary}[theorem]{Corollary}
\theoremstyle{definition}

\theoremstyle{remark}
\newtheorem{remark}[theorem]{Remark}

\begin{document}
\title[Maximal curves with cyclic automorphism groups]
{New families of maximal curves\\
not covered by the Hermitian curve\\
and their cyclic automorphism groups}
\author{Xubin Hu}
\address{School of Mathematical Sciences, University of Science and Technology of China,
Hefei, Anhui 230026, China}
\email{xbh5168@mail.ustc.edu.cn}
\author{Liming Ma}
\address{School of Mathematical Sciences, University of Science and Technology of China,
Hefei, Anhui 230026, China}
\email{lmma20@ustc.edu.cn}
\author{Yipeng Wang}
\address{School of Mathematical Sciences, University of Science and Technology of China,
Hefei, Anhui 230026, China}
\email{jisunxingfu@mail.ustc.edu.cn}
\date{November 7, 2026}
\hypersetup{pdftitle={New families of maximal curves not covered by the Hermitian curve and their cyclic automorphism groups},pdfauthor={Xubin Hu, Liming Ma, Yipeng Wang}}

\begin{abstract}
For every prime power $q>2$, integer $n\ge2$ and $1\le r\le2n-1$,
we construct an $\F_{q^{2n}}$-maximal curve $T_{q,n,r}$ with
Kummer equation
\[
 z^{q^n+1}=\eta^{-R}y(1-y^{q-1})^R,
 \qquad R=\frac{q^r-1}{q-1},
\]
where $\eta\in\F_{q^{2n}}^*$ satisfies $\eta^{q^n}=-\eta$.
For fixed $q,n$, these are cyclic covers of the projective line
over a single field, indexed by the integer $r$.
Writing $\delta=\gcd(q^n+1,R)$, we prove that the genus is
$(q-1)(q^n+1-\delta)/2$ and that $T_{q,n,r}$ is covered by
the Hermitian curve $\HH_{q^n}$ over an algebraic closure if and only if
$\min(r,2n-r)$ divides $n$.
For $1<r<n$ and $\delta=1$, the full automorphism group is
cyclic of order $(q-1)(q^n+1)$, and distinct values of $r$ give
curves that are nonisomorphic over an algebraic closure despite
having the same genus and full automorphism group.
We also give explicit equations and genera for all their Galois
subfields. When $n\ge4$ is a power of two, this produces families
over a fixed field whose members are not covered by $\HH_{q^n}$.
\end{abstract}
\maketitle
\enlargethispage{6pt}

\section{Introduction}
\label{sec:introduction}

Let $Q$ be a power of a prime. A nonsingular, geometrically irreducible projective curve $X$ over
$\F_{Q^2}$ is called maximal if its number of rational points attains
the Hasse--Weil upper bound
$\#X(\F_{Q^2})=Q^2+1+2g(X)Q$, where $g(X)$ is the genus of the curve $X$. The Hermitian curve $\HH_Q$, with affine equation
$b^Q+b=a^{Q+1}$, is maximal and has genus $Q(Q-1)/2$.
By R\"uck and Stichtenoth \cite{RS1994}, it is the unique
$\F_{Q^2}$-maximal curve of this genus, up to $\F_{Q^2}$-isomorphism.
A theorem of Serre states that a curve covered over $\F_{Q^2}$ by
an $\F_{Q^2}$-maximal curve is also maximal; see
\cite[Proposition~6]{Lachaud1987}.
Thus subfields of the Hermitian function field provide a class of
maximal function fields whose genera and defining equations have
been studied in \cite{CKT2000,GSX2000}.

Constructing maximal curves which are not covered by the Hermitian
curve is more difficult: one must attain the Hasse--Weil bound
and rule out every covering from the Hermitian curve, including
non-Galois coverings.
Giulietti and Korchm\'aros constructed maximal curves which are not
covered by the Hermitian curve over the same finite field
\cite{GK2009}. Garcia, G\"uneri and Stichtenoth \cite{GGS2010},
and later Beelen and Montanucci \cite{BM2018}, constructed two
families extending the GK curve to $\F_{q^{2n}}$ for odd $n$.
Ma and Wang \cite{MW2026} gave a Kummer model which also defines
maximal curves for even $n$, and proved that its members with
$q>2$ and $n\ge3$ are not covered by $\HH_{q^n}$ over
$\F_{q^{2n}}$.
They also proved the corresponding non-covering result for the GGS
curves when $q>2$ and $n\ge3$ is odd
\cite[Theorem~3.5]{MW2026}.

In this paper we construct new families of maximal curves not
covered by the Hermitian curve. The construction uses a Kummer
equation with an integer parameter $r$ that can vary while the field
$\F_{q^{2n}}$ remains fixed.
We prove maximality throughout the full parameter range and
determine exactly which members are covered by the Hermitian curve.
For an explicit range of parameters, the full automorphism groups
are cyclic, and distinct values of $r$ give nonisomorphic curves
with the same genus and full automorphism group.
We also give explicit equations and genera for all their Galois subfields.

\subsection{Our results}

Throughout, a curve means a
nonsingular geometrically integral projective variety of dimension one.

Let $q>2$ be a prime power, let $n\ge2$, and let $1\le r\le2n-1$.
Put
\begin{equation}
\label{eq:parameters}
 Q=q^n,\qquad N=Q+1,\qquad R=\frac{q^r-1}{q-1},
 \qquad k=\F_{q^{2n}},
\end{equation}
and write $\delta=\gcd(N,R)$.
Choose $\eta\in k^*$ with $\eta^Q=-\eta$.
We denote by $T_{q,n,r}$ the nonsingular projective curve with
function field
\begin{equation}
\label{eq:model}
 k(y,z),\qquad z^N=\eta^{-R}y(1-y^{q-1})^R.
\end{equation}
The extension $k(T_{q,n,r})/k(y)$ is cyclic of degree $N$,
with precisely the $q+1$ rational places $y\in\F_q\cup\{\infty\}$
ramified. Thus the parameter $r$ changes the exponents while
leaving the set of ramified places fixed.
Passing to $t=y^{q-1}$ gives a cyclic extension of degree
$(q-1)N$ with only three ramified places, as in
\eqref{eq:three-branch}; this model will be used to determine
the automorphism groups and distinguish the isomorphism classes.

The $k$-isomorphism class is independent of $\eta$, and
$T_{q,n,r}\cong_kT_{q,n,2n-r}$ by
Proposition~\ref{prop:parameter-symmetry}, so it suffices to
consider $1\le r\le n$.
Throughout, $\Aut(T)$ denotes the full automorphism group over
an algebraic closure $\overline{k}$ of $k$; for a maximal curve
of genus at least two, all these automorphisms are defined over
the field of maximality \cite[Theorem~3.10]{GSY2015}.

We first establish maximality and compute the genus.
\begin{theorem}
\label{thm:A}
The curve $T_{q,n,r}$ is maximal over $k$, with
\begin{equation}
\label{eq:intro-genus-count}
 g(T_{q,n,r})=\frac{(q-1)(N-\delta)}2,\qquad
 \#T_{q,n,r}(k)=Q^2+1+(q-1)Q(N-\delta).
\end{equation}
\end{theorem}

A covering over a field $\mathbb F$ means a nonconstant
$\mathbb F$-morphism between curves. The Hermitian coverings admit the following characterization.

\begin{theorem}
\label{thm:non-cover-by-Herm}
The curve $T_{q,n,r}$ is covered by $\HH_Q$ over $\overline{k}$
if and only if $\min(r,2n-r)$ divides $n$.
Whenever such a covering exists, one can choose it Galois and
defined over $k$.
\end{theorem}

The construction includes Hermitian quotients and two previously studied families.
When $r\mid n$, Proposition~\ref{prop:explicit-cover} realizes
$T_{q,n,r}$ as a $k$-Galois quotient of $\HH_Q$ of degree
$q^{n-r}R$.
For odd $n\ge3$, the member $T_{q,n,2}$ is the curve
$\mathcal X_n$ whose maximality was proved by Abd\'on, Bezerra
and Quoos \cite{ABQ2009} (for prime $q$, see also Garcia and
Tafazolian \cite{GT2008}). It is a degree-$q$ quotient of the
Garcia--G\"uneri--Stichtenoth curve
\cite[Proposition~2.2]{GGS2010}; see Remark~\ref{rem:GS}.
For the same odd $n$, the function field of $T_{q,n,n-1}$ is
the subfield $k(\rho,w)$ of the Beelen--Montanucci function field
explicitly described and proved maximal in \cite{BM2018}; see
Remark~\ref{rem:endpoint}.

Theorem~\ref{thm:non-cover-by-Herm} shows that these two degree-$q$ quotients
are themselves not covered by $\HH_Q$, even over $\overline{k}$.
For the following curves $\mathcal X_n$, this settles the covering question
recorded in \cite[Introduction]{DM2012}. We state this consequence
separately.

\begin{corollary}
\label{cor:Xn}
Let $q>2$ and let $n\ge3$ be odd. Let $\mathcal X_n$ be the
nonsingular projective model of the plane curve
\[
 z^{(q^n+1)/(q+1)}=x^{q^2}-x.
\]
Then $\mathcal X_n$ is not covered by $\HH_{q^n}$ over
$\overline{\F}_{q^{2n}}$.
\end{corollary}
\begin{proof}
Using the isomorphism $\mathcal X_n\cong_k T_{q,n,2}$ established
in Remark~\ref{rem:GS}, the assertion follows from
Theorem~\ref{thm:non-cover-by-Herm} and $2\nmid n$.
\end{proof}

We also study coverings of $T_{q,n,r}$ by the curves in
\cite{MW2026} for fixed $q,n$; the precise parameter ranges and
covering criterion are given in Theorem~\ref{thm:main}.

Another question is whether different parameters give
isomorphic curves. Genus and automorphism group alone need not
distinguish maximal function fields: Beelen, Montanucci, Niemann
and Quoos \cite{BMNQ2025} obtained nonisomorphic Hermitian
Galois subfields with the same genus and, in many cases, the
same automorphism group.
Beelen, Drue, Montanucci and Zini \cite{BDMZ2026} distinguished
related BM and GGS subfields by their automorphism groups.
In our construction, varying $r$ gives curves with the same
genus and full automorphism group which are nevertheless
nonisomorphic over $\overline{k}$, under the condition
\begin{equation}
\label{eq:range}
 1<r<n,\qquad \gcd\left(q^n+1,\frac{q^r-1}{q-1}\right)=1.
\end{equation}
In this range the genus is $(q-1)Q/2$ and
$\#T_{q,n,r}(k)=qQ^2+1$.
Every odd $r$ with $1<r<n$ satisfies \eqref{eq:range}.

\begin{theorem}
\label{thm:B}
Assume \eqref{eq:range}, and put $M=(q-1)(q^n+1)$.
The automorphisms of $T_{q,n,r}$ are exactly
\begin{equation}
\label{eq:intro-aut}
 \sigma_\lambda(y)=\lambda^{q^n+1}y,\qquad
 \sigma_\lambda(z)=\lambda z,\qquad \lambda^M=1.
\end{equation}
Thus $\Aut(T_{q,n,r})\cong C_M$.
For fixed $q,n$ and parameters $r,r'$ satisfying \eqref{eq:range},
\begin{equation}
\label{eq:intro-isomorphism}
 T_{q,n,r}\cong_{\overline{k}}T_{q,n,r'}
 \quad\Longleftrightarrow\quad r=r'.
\end{equation}
\end{theorem}

For $n$ a power of two,
Theorems~\ref{thm:A},~\ref{thm:non-cover-by-Herm} and~\ref{thm:B}
yield the following families of pairwise nonisomorphic maximal curves
over a single field.

\begin{corollary}
\label{cor:C}
Fix a prime power $q>2$. Let $n=2^m$, where $m\ge2$.
Let
\begin{equation}
\label{eq:power-two-family}
 \mathcal R=
 \begin{cases}
 \{3,5,\ldots,2^m-1\},&q\ \text{odd},\\
 \{r:1<r<2^m,\ r\nmid2^m\},&q\ \text{even}.
 \end{cases}
\end{equation}
The curves $T_{q,n,r}$ with $r\in\mathcal R$ are maximal over
$\F_{q^{2n}}$ and pairwise nonisomorphic over
$\overline{\F}_{q^{2n}}$. They all have genus $(q-1)q^n/2$ and full
automorphism group $C_{(q-1)(q^n+1)}$. None is covered by
$\HH_{q^n}$ over $\overline{\F}_{q^{2n}}$.
There are $2^{m-1}-1$ such curves when $q$ is odd and
$2^m-m-1$ when $q$ is even.
\end{corollary}

\subsection{Methods}

We prove maximality by constructing two disjoint sets of rational
places of $k(y)$ that split completely in $k(T_{q,n,r})$.
Together with the ramified places, they give the Hasse--Weil upper bound.
For non-covering by Hermitian curves, we combine the polynomial representation in
\cite[Lemma~3.1]{MW2026} with a local restriction on vanishing orders
on the Hermitian curve.
For $1<r<n$, a hypothetical covering produces a polynomial
$H(a,b)$ on $\HH_Q$ with $\deg H\le q^{n-r}$.
Writing $A=q^r$, comparison of coefficients forces nonzero terms
at the successive powers $b,b^A,b^{A^2},\ldots$.
The Hermitian equation can terminate this recurrence only at
$b^Q$; otherwise the recurrence exceeds the degree bound for $H$.
This requires $A^j=Q$ for some $j$, and hence $r\mid n$.

To determine the automorphism group, we compute the normalizer of the
Kummer group of order $N$ and the stabilizers of its fixed points.
The Hurwitz genus formula then gives the full group. Its fixed field is
rational, and exactly three places of this field ramify in
$k(T_{q,n,r})$.
Comparing their Kummer exponents proves
\eqref{eq:intro-isomorphism}.
The same three-place model gives all Galois subfields and their genera.

\subsection{Organization}
Section~\ref{sec:family} gives the models, computes the genus,
proves maximality and constructs the Hermitian coverings.
Section~\ref{sec:non_cover_by_X} studies coverings by the curve
$\X_{q,n}$ of \cite{MW2026}.
Section~\ref{sec:noncover} proves the covering criterion in
Theorem~\ref{thm:non-cover-by-Herm}.
Section~\ref{sec:automorphisms} determines the full automorphism
group, isomorphism classes and Galois subfields under \eqref{eq:range}.

\section{A family of maximal curves}
\label{sec:family}

We retain the notation of the introduction and write $p$ for the
characteristic of $k$. We use the notation and terminology of
\cite{Stichtenoth2009} for algebraic function fields and
\cite{Har77, HKT2008} for algebraic curves.
All isomorphisms and morphisms are over $k$ unless another field
is specified. A place is rational if its residue field is $k$; places over
$\overline{k}$ are identified with geometric points.
For a nonzero function $f$, write $(f)_0$ and $(f)_\infty$ for
its zero and pole divisors.

\subsection{Construction and genus}

The right-hand side of \eqref{eq:model} has a simple zero at $y=0$.
Eisenstein's criterion at $y=0$, applied over $\overline{k}$,
shows that the defining polynomial is absolutely irreducible and
$[k(y,z):k(y)]=N$. Since $p\nmid N$ and $N\mid Q^2-1$, this
extension is Kummer and cyclic of degree $N$.

If $\eta'\in k^*$ is another choice with $\eta'^Q=-\eta'$, then
$a=\eta'/\eta\in\F_Q^*$. Surjectivity of the norm
$k^*\to\F_Q^*$ gives $\lambda\in k^*$ with $\lambda^N=a^{-R}$.
Replacing $z$ by $\lambda z$ gives the equation with coefficient
$\eta'^{-R}$. Thus the $k$-isomorphism class is independent of $\eta$.

The ramification depends on $\delta=\gcd(N,R)$, whose value is
given by the following lemma.

\begin{lemma}
\label{lem:arithmetic}
Let $h=\gcd(n,r)$. Then
\begin{equation}
\label{eq:delta}
 \delta=
 \begin{cases}
 q^h+1,&r/h\ \text{even},\\
 2,&r/h\ \text{odd},\ q\ \text{odd},\ r\ \text{even},\\
 1,&\text{otherwise}.
 \end{cases}
\end{equation}
If $1<r<n$, then $R>1$ and
\begin{equation}
\label{eq:arithmetic}
 R+q^r=\frac{q^{r+1}-1}{q-1}<N.
\end{equation}
\end{lemma}
\begin{proof}
The gcd of $q^n+1$ and $q^r-1$ is $q^h+1$ when $r/h$ is even,
and is at most two otherwise.
In the first case, writing
\[
 R=\frac{q^h-1}{q-1}
   \bigl(1+q^h+\cdots+q^{h(r/h-1)}\bigr)
\]
shows that $q^h+1$ divides $R$.
In the other case, $N$ is odd when $q$ is even; when $q$ is
odd, $N$ is even and $R\equiv r\pmod2$.
This gives \eqref{eq:delta}.
For $1<r<n$, we have $R>1$ and
$R+q^r\le(Q-1)/(q-1)<N$.
\end{proof}

\begin{lemma}
\label{lem:ramification}
The places $y=0,\infty$ are totally ramified in $k(y,z)/k(y)$;
write $P_0,P_\infty$ for their extensions.
For each $\alpha\in\F_q^*$, there are $\delta$ rational places
$P_{\alpha,1},\ldots,P_{\alpha,\delta}$ above $y=\alpha$,
each with ramification index $N/\delta$.
There is no ramification elsewhere, and
\begin{equation}
\label{eq:divisors}
 (y)=N(P_0-P_\infty),\qquad
 (z)=P_0+\frac R\delta
       \sum_{\alpha\in\F_q^*}\sum_{j=1}^{\delta}P_{\alpha,j}
       -q^rP_\infty.
\end{equation}
Moreover,
\begin{equation}
\label{eq:genus}
 2g(T_{q,n,r})-2=(q-1)(N-\delta)-2.
\end{equation}
When $\delta=1$, write $P_\alpha=P_{\alpha,1}$.
\end{lemma}
\begin{proof}
The valuations of the right-hand side of \eqref{eq:model} at
$0,\infty,\alpha\in\F_q^*$ are $1,-q^r,R$.
Proposition~3.7.3 in \cite{Stichtenoth2009} gives the stated
ramification indices.
The places $P_0,P_\infty$ are rational by total ramification.
The function
\[
 v=\frac{\eta^{R/\delta}z^{N/\delta}}
          {(1-y^{q-1})^{R/\delta}}
 \qquad\text{satisfies}\qquad v^\delta=y.
\]
Every $\alpha\in\F_q^*$ has $\delta$ distinct $\delta$-th roots in
$k$, since $\delta\mid N$ and $q-1\mid(Q^2-1)/\delta$.
Thus $y=\alpha$ splits in $k(v)/k(y)$; each resulting place is
totally ramified in $k(y,z)/k(v)$, so its residue field remains $k$.
The valuations give \eqref{eq:divisors}.
The Hurwitz genus formula gives
\[
 2g-2=-2N+2(N-1)+(q-1)(N-\delta),
\]
which is \eqref{eq:genus}.
\end{proof}

Put $t=y^{q-1}$ and $M=(q-1)N$. Equation~\eqref{eq:model} gives
\begin{equation}
\label{eq:three-branch}
 z^M=\eta^{1-q^r}t(1-t)^{q^r-1},\qquad
 y=\frac{\eta^Rz^N}{(1-t)^R}.
\end{equation}
Thus $k(y,z)=k(t,z)$, and $k(t,z)/k(t)$ is cyclic of degree $M$.
The only places of $k(t)$ ramified in this extension are
$t=0,1,\infty$, with indices $M,N/\delta,M$, respectively, because
\begin{equation}
\label{eq:gcd-M}
 \gcd(M,q^r-1)=(q-1)\delta.
\end{equation}

\begin{proposition}
\label{prop:parameter-symmetry}
There is a $k$-isomorphism
$T_{q,n,r}\cong T_{q,n,2n-r}$.
\end{proposition}
\begin{proof}
Let $r'=2n-r$ and $D_0=(Q-1)/(q-1)$.
In \eqref{eq:three-branch}, set
\[
 t'=\frac1t,\qquad z'=\frac{(t-1)^{D_0}}{z^{q^{r'}}}.
\]
Since $q^rq^{r'}=Q^2=1+D_0M$ and $\eta^{Q^2}=\eta$, substitution gives
\[
 (z')^M=\eta^{1-q^{r'}}t'(1-t')^{q^{r'}-1}.
\]
As $\gcd(q^{r'},M)=1$, the new functions generate $k(t,z)$.
Their equation is \eqref{eq:three-branch} with parameter $r'$, so
this field is $k(T_{q,n,r'})$.
\end{proof}

\subsection{Maximality}

\begin{proposition}
\label{prop:maximality}
For every $1\le r\le2n-1$, the curve $T_{q,n,r}$ is maximal
over $k$, and
\begin{equation}
\label{eq:point-count}
 \#T_{q,n,r}(k)=Q^2+1+(q-1)Q(N-\delta).
\end{equation}
In particular, $\#T_{q,n,r}(k)=qQ^2+1$ when $\delta=1$.
\end{proposition}
\begin{proof}
Put $f(Y)=\eta^{-R}Y(1-Y^{q-1})^R$ and $A=q^r$.
The norm map $k^*\to\F_Q^*$, $a\mapsto a^N$, is surjective
with kernel of order $N$; see \cite[Chapter~2]{LN1997}.
Thus, for $\beta\in k\setminus\F_q$, the place $y=\beta$ splits
completely in $k(T_{q,n,r})/k(y)$ if and only if
$f(\beta)\in\F_Q^*$.
We construct two disjoint sets of such values $\beta$ and count
the resulting rational places.

Let $D=\F_{q^{\gcd(2n,r)}}$, and put
\[
 U=\{u\in k:u^Q+u=1\},\qquad
 L=\{\beta\in k:\beta^Q=(-1)^R\beta\}.
\]
Both sets have $Q$ elements. The first construction uses the
trace equation defining $U$: for $u\in U\setminus D$ and
$\alpha\in\F_q^*$, we have
\begin{equation}
\label{eq:split-values}
 f(\alpha u^R)=\alpha\{\eta^{-1}(u-u^A)\}^R\in\F_Q^*.
\end{equation}
Indeed, $u^A\ne u$ and $(u-u^A)^Q=-(u-u^A)$, so division
by $\eta$ places $\eta^{-1}(u-u^A)$ in $\F_Q^*$.
Since $f$ vanishes precisely on $\F_q$, each $\alpha u^R$
lies outside $\F_q$. To count these values, we check that
different pairs $(\alpha,u)$ give different elements.
Indeed, if $\alpha u^R=\alpha'v^R$,
then $\lambda=u/v$ satisfies $\lambda^{A-1}=1$, so
$\lambda\in D$. The equations defining $U$ give
\[
 (\lambda-\lambda^Q)v=1-\lambda^Q.
\]
If $\lambda=\lambda^Q$, then $\lambda=1$; otherwise $v\in D$,
a contradiction. Hence $u=v$ and $\alpha=\alpha'$.

The second construction takes $\beta$ directly from
$L\setminus\F_q$. Here $\beta^{q-1}\in\F_Q$ and
$f(\beta)^Q=f(\beta)$, so these places also split completely.

The two sets of places are disjoint. For if $\alpha u^R\in L$, then
$\lambda=-u^Q/u=(u-1)/u$ satisfies $\lambda^R=1$.
Thus $\lambda\in D$, and $u=1/(1-\lambda)\in D$, a contradiction.

Write $\tau=|U\cap D|$, and let $\epsilon=1$ if $p=2$ or $r$
is even, and $\epsilon=0$ otherwise.
Then $|L\setminus\F_q|=Q-1-\epsilon(q-1)$.
To compute $\tau$, put $h=\gcd(n,r)$.
If $r/h$ is even, then $D=\F_{q^{2h}}$ and $n/h$ is odd,
so $u^Q+u=1$ on $D$ is the trace equation $u^{q^h}+u=1$;
hence $\tau=q^h$.
If $r/h$ is odd, then $D=\F_{q^h}\subseteq\F_Q$, and the
equation becomes $2u=1$; hence $\tau=1$ for odd $p$ and $\tau=0$
for $p=2$.
Lemma~\ref{lem:arithmetic} now gives $\tau+\epsilon=\delta$.

There are $(q-1)(Q-\tau)$ values from the first construction
and $Q-1-\epsilon(q-1)$ from the second, each yielding $N$
rational places. Lemma~\ref{lem:ramification} supplies a further
$2+(q-1)\delta$ places above $\F_q\cup\{\infty\}$.
Together these contributions give
\begin{align*}
 \#T_{q,n,r}(k)
 &\ge N(q-1)(Q-\tau)
       +N\bigl(Q-1-\epsilon(q-1)\bigr)+2+(q-1)\delta\\
 &=Q^2+1+(q-1)Q(N-\delta).
\end{align*}
By \eqref{eq:genus}, the last expression is the Hasse--Weil
upper bound, so equality holds and the curve is maximal.
\end{proof}

Together with Lemma~\ref{lem:ramification},
Proposition~\ref{prop:maximality} proves Theorem~\ref{thm:A}.

\subsection{Hermitian subcovers and known members}

\begin{proposition}
\label{prop:explicit-cover}
If $r\mid n$, there is a $k$-Galois covering
$\HH_Q\to T_{q,n,r}$ of degree $q^{n-r}R$.
\end{proposition}
\begin{proof}
On the Hermitian model
\[
 V^N=\eta^{-1}(X-X^Q),
\]
obtained from $b^Q+b=a^N$ by $X=\eta b$, $V=a$, put
\begin{equation}
\label{eq:explicit-Hermitian-cover}
 U=\sum_{j=0}^{n/r-1}X^{q^{rj}},\qquad
 y=U^R,\qquad z=V^R.
\end{equation}
Since $U-U^{q^r}=X-X^Q$, these functions satisfy
\eqref{eq:model}.
The transformations
\[
 (X,V)\longmapsto(\zeta^N X+a,\zeta V),
 \qquad \zeta^R=1,\quad
 a\in\ker\Tr_{\F_Q/\F_{q^r}},
\]
form a group of order $q^{n-r}R$ fixing $y,z$.
The trace step has degree $q^{n-r}$.
Moreover,
$U=\eta V^N/(1-y^{q-1})$, so the remaining field extension
is generated by $V$ with $V^R=z$ and has degree at most $R$.
The $q^{n-r}R$ automorphisms fixing $y,z$ force equality.
This is the required $k$-Galois covering.
\end{proof}

\begin{remark}
\label{rem:GS}
For odd $n\ge3$, the curve $T_{q,n,2}$ is $k$-isomorphic
to the curve
\[
 z^{(Q+1)/(q+1)}=x^{q^2}-x.
\]
Indeed, choose $\eta\in\F_{q^2}^*$ with $\eta^q=-\eta$,
put $m=(Q+1)/(q+1)$ and
$x=-z^m/(1-y^{q-1})$.
Then $y=\eta^{q+1}x^{q+1}$ and $z^m=x^{q^2}-x$.
This is the curve $\mathcal X_n$ in Corollary~\ref{cor:Xn}.
The GGS curve is obtained by adjoining $v$ with
$v^q+v=x^{q+1}$; see \cite[equation~(2.1)]{GGS2010}.
Forgetting $v$ gives a Galois covering of degree $q$, with group
$v\mapsto v+a$, $a^q+a=0$.
Corollary~\ref{cor:Xn} excludes every covering
$\HH_Q\to T_{q,n,2}$ over $\overline{k}$.
For $n=3$, it was proved that no Galois covering
$\HH_{q^3}\to T_{q,3,2}$ exists over $\F_{q^6}$,
for $q=3$ in \cite{GS2006} and for $q>3$ in \cite{GMZ2016}.

This example also shows that the hypothesis $\delta=1$ in
Theorem~\ref{thm:B} cannot be omitted.
Here $\delta=q+1$, and the translations
$(x,z)\mapsto(x+\beta,z)$, $\beta\in\F_{q^2}$,
form an automorphism group of order $q^2$.
Since $p\nmid(q-1)(Q+1)$, the full group strictly contains
the cyclic subgroup in \eqref{eq:intro-aut}.
We do not determine the full groups for general $\delta>1$.
\end{remark}

\begin{remark}
\label{rem:endpoint}
For $n\ge3$, the curve $T_{q,n,n-1}$ is the quotient of the
curve in \cite{MW2026} by its group of $q$ translations.
For odd $n$, the upper curve is the Beelen--Montanucci curve;
for even $n$, it belongs to the family constructed in \cite{MW2026}.
To describe the covering, put $r=n-1$, $A(u)=u^q-u$, and choose
$\theta\in\F_{q^2}^*$ with $\theta^q=-\theta$.
Then $c=\eta^{-R}\theta$ satisfies $c^Q=-c$.
Indeed, in odd characteristic $R\equiv n-1\pmod2$ and
$\theta^Q=(-1)^n\theta$, so $c^Q=(-1)^{R+n}c=-c$;
in characteristic two the equality is immediate.
The curve in \cite{MW2026} has the model
\[
 Z^N=cA(u)\bigl(1+A(u)^{q-1}\bigr)^R.
\]
Since $\theta^{q-1}=-1$, the functions
$y=\theta A(u)$ and $z=Z$ satisfy \eqref{eq:model}.
The translations $u\mapsto u+a$, $Z\mapsto Z$, with
$a\in\F_q$, fix $y,z$. Since $u^q-u=y/\theta$, their fixed
field is $k(y,z)$ and the covering has degree $q$.
For $n=3$, this covering and the GGS covering in
Remark~\ref{rem:GS} both have the GK curve as source.

For odd $n$, this quotient was already described in \cite{BM2018},
where its maximality and genus were proved. It is the fixed field
of the group of order $q$ considered in
\cite[Lemma~3.3]{BDMZ2026}.
To identify this subfield explicitly, assume $n$ is odd and put
\[
 m=\frac{Q+1}{q+1},\qquad
 s=\frac{R}{q+1}=\frac{q^{n-1}-1}{q^2-1},
 \qquad q(q-1)s=m-1.
\]
By Lemma~\ref{lem:ramification}, the function
$v=\eta^s z^m/(1-y^{q-1})^s$ satisfies $v^{q+1}=y$, and
\[
 k(T_{q,n,n-1})=k(v,z),\qquad
 z^m=\eta^{-s}v(1-v^{q^2-1})^s.
\]
Choose $\xi\in k^*$ with $\xi^m=\eta^{1-m}$.
Such a choice is possible: $m$ is odd, so
$\eta^{1-m}\in\F_Q^*$, and $m\mid Q+1$ implies
$\F_Q^*\subseteq(k^*)^m$.
Then
\[
 \rho=v^{-1},\qquad
 w=\frac{\xi(1-v^{q^2-1})}{z^{q(q-1)}}
 \quad\text{satisfy}\quad
 w^m=\frac{\rho^{q^2-1}-1}{\rho^{q-1}}.
\]
Conversely, $v=\rho^{-1}$ and
$z=(\eta\xi)^{-s}\rho^{-1}w^s$, so $k(v,z)=k(\rho,w)$.
This is precisely the subfield in \cite{BM2018}.

For either parity of $n$, Theorem~\ref{thm:non-cover-by-Herm} proves that the
quotient $T_{q,n,n-1}$ itself is not covered by $\HH_Q$ over
$\overline{k}$. This strengthens the non-covering result for
the upper curve in \cite{MW2026}.
\end{remark}

\section{Coverings by the curves \texorpdfstring{$\mathcal{X}_{q,n}$}{X\_\{q,n\}}} \label{sec:non_cover_by_X}
In this section we study coverings of $T_{q,n,r}$ by the curve
$\X_{q,n}$ of \cite{MW2026} for the same parameters $q,n$.
Choose $c\in k^*$ with $c^Q=-c$. The curve $\X_{q,n}$ has the model
\begin{equation}
     Z^N =cA(u)B(u)^{R_0},\label{eq:X}
\end{equation}
where $A(u)=u^q-u$, $B(u)=1+A(u)^{q-1}$ and
$R_0=(q^{n-1}-1)/(q-1)$.

The main result is as follows. 

\Needspace{10\baselineskip}
\begin{theorem}\label{thm:main}
Let $q>2$ and $1<r<n$.

\begin{enumerate}[label=(\arabic*)]
\item If $n\ge6$ is even, then $T_{q,n,r}$ is covered by $\X_{q,n}$
over $\overline{k}$ if and only if $r=n-1$. Every separable covering in
this case is Galois of degree $q$ and is defined over $k$. Its
Galois group is conjugate in $\Aut(\X_{q,n})$ to the group of translations
$(u,Z)\mapsto(u+a,Z)$, $a\in\F_q$.
\item If $n\ge5$ is odd and $\delta=1$, then $T_{q,n,r}$ is not
covered by $\X_{q,n}$ over $\overline{k}$.
\end{enumerate}
\end{theorem}

We use the following descent property for separable morphisms
between maximal curves.

\begin{proposition}
\label{prop:descent}
Let $X,Y$ be maximal curves over $k=\F_{Q^2}$, with $g(Y)\ge2$.
Every nonconstant separable $\overline{k}$-morphism $f:X\to Y$
is defined over $k$.
\end{proposition}
\begin{proof}
Put $d=\deg f$. The Hurwitz genus formula gives $d\le g(X)-1$,
whereas $\#X(k)=Q^2+1+2Qg(X)>2d$.
Thus there are three points $P_1,P_2,P_3\in X(k)$ with
distinct images under $f$.
The Frobenius divisor relation
\cite[Proposition~3.7]{GSY2015} gives
$(Q+1)P_1\sim(Q+1)P_2\sim(Q+1)P_3$.
Pushforward preserves linear equivalence, since the pushforward
of a principal divisor is the divisor of the norm of its function.
Thus the $(Q+1)$-multiples of the three image points are linearly
equivalent, and \cite[Proposition~3.8]{GSY2015} shows that
these points belong to $Y(k)$.

Write $\Fr_X,\Fr_Y$ for the $Q^2$-power Frobenius maps.
For any geometric point $P$ of $X$, push forward
$QP+\Fr_X(P)\sim(Q+1)P_1$ and compare with the
Frobenius divisor relation on $Y$. This gives
\[
 Qf(P)+f(\Fr_X(P))\sim(Q+1)f(P_1)
 \sim Qf(P)+\Fr_Y(f(P)).
\]
Two linearly equivalent point divisors on a curve of positive
genus are equal. Hence $f$ commutes with the $Q^2$-power
Frobenius and is defined over $k$.
\end{proof}

Write $\X=\X_{q,n}$ and $T=T_{q,n,r}$.
By \cite{MW2026} and Theorem~\ref{thm:A}, both curves are maximal
over $k$, with genera
\begin{equation}
    2g(\X)=(q^2-1)Q-\epsilon_n q^2(q-1),\qquad
2g(T)=(q-1)(N-\delta), \label{eq:genera}
\end{equation}
where $\epsilon_n = 0$ for even $n$ and $\epsilon_n = 1$ for odd $n$. 

The existence of a covering over $\overline{k}$ reduces to that
of a separable covering $\pi:\X\to T$ over $k$.
Indeed, every Frobenius twist of $\X$ has model \eqref{eq:X}
with another nonzero coefficient, so scaling $Z$ gives a
$\overline{k}$-isomorphism with $\X$.
Factoring through a power of the relative Frobenius therefore
gives a separable covering $\X\to T$ over $\overline{k}$.
Proposition~\ref{prop:descent} makes it defined over $k$,
since $g(T)>1$ by \eqref{eq:genera}.
In the proof of Theorem~\ref{thm:main}, we therefore assume that
$\pi:\X\to T$ is a nonconstant separable covering over $k$.

Put $F=\overline{k}(\X_{q,n})$, $b=q^2-q$, and
\begin{equation}
K=q(q^2-1),\qquad w=B(u)/Z^b .
\label{eq:w}
\end{equation}
We use degrees of rational functions over $\overline{k}$, so that
$\deg f=[F:\overline{k}(f)]$ for $f\in F$ nonconstant.

\begin{lemma}\label{lem:rational}
The extension $F/\overline{k}(w)$ is Galois with group
$G_0\cong \SL_2(q)$ and degree $K$.
\end{lemma}
\begin{proof}
Each $M=\left(\begin{smallmatrix}\alpha&\beta\\ \gamma&d_0\end{smallmatrix}\right)
\in\SL_2(q)$ acts on the curve by the substitution
\[
(u,Z)\longmapsto
\left(Mu,\frac{Z}{\gamma u+d_0}\right),\qquad
Mu=\frac{\alpha u+\beta}{\gamma u+d_0}.
\]
Indeed, $A(Mu)=A(u)/(\gamma u+d_0)^{q+1}$ and
$B(Mu)=B(u)/(\gamma u+d_0)^b$. Since $q+1+bR_0=N$, the substitution
preserves \eqref{eq:X} and fixes $w$. This action is faithful, since only scalar matrices
fix $u$, and $-I$ sends $Z$ to $-Z$ in odd characteristic.

To compute the extension degree, observe that
\[
 w^N=c^{-b}\frac{B(u)^{q+1}}{A(u)^b}.
\]
The right side has degree $K$ as a rational function of $u$ because
$A,B$ are relatively prime and the numerator has degree $b(q+1)=K$.
Since $[F:\overline{k}(u)]=N=[\overline{k}(w):\overline{k}(w^N)]$,
the tower law gives $[F:\overline{k}(w)]=K= |G_0|$.
\end{proof}

\begin{lemma}\label{lem:small}
Let $n\ge4$ be even. If $f\in F$ is nonconstant and
$2q\deg f\le Q-1$, then $f\in\overline{k}(w)$. In particular,
$K\mid\deg f$.
\end{lemma}
\begin{proof}
Suppose $f\notin\overline{k}(w)$ and write $D=\deg f$.
By Lemma~\ref{lem:rational}, $L=\overline{k}(f,w)=F^H$
for a proper subgroup $H<G_0$. Put $e=|H|$.
The two rational subfields generating $L$ have indices $D/e$ and
$K/e$. The Castelnuovo--Severi inequality
\cite[Theorem~3.11.3]{Stichtenoth2009} gives
\[
g(L)\le(D/e-1)(K/e-1)<DK/e^2.
\]
Let $S$ be the image of $H$ in $\PGL_2(q)$, let $\nu(S)$ count its
elements of order $p$, and put $a=q^2-1-\nu(S)$.
The scalar subgroup of $H$ has order $c_H\le2$, so the quotient
genus formula \cite[Theorem~C(2)]{MW2026} gives $2g(L)=a(N-c_H)/e$.
Consequently
\begin{equation}
e a(Q-1)<2DK.
\label{eq:bound}
\end{equation}

We claim that $a\ge1$. The projective image of $G_0$ has $q^2-1$
elements of order $p$, namely $q-1$ for each of the $q+1$ possible
fixed points in $\mathbb P^1(\F_q)$. If all belong to $S$, then
$H$ contains their unipotent lifts. Indeed, in odd characteristic
a lift in $H$ of a projective unipotent is $U$ or $-U$, where
$U$ is unipotent, and $(-U)^{p+1}=U$. In characteristic two the
lift is unique. Thus $H$ would contain the upper and lower unipotent
groups, which generate $G_0$. This contradicts $H<G_0$.

Now $\nu(S)\le e-1$, so $a+e\ge q^2$. Since $a,e\ge1$, we obtain
$ea\ge e+a-1\ge q^2-1$. Substituting in \eqref{eq:bound} gives
$Q-1<2qD$. This proves the first assertion. The second follows
from $[F:\overline{k}(w)]=K$.
\end{proof}

\begin{proof}[Proof of Theorem~\ref{thm:main}\textup{(1)}]
Write $d=\deg\pi$. Since
$(q-1)\delta\le(q-1)R\le Q/q-1$, the Hurwitz genus formula
and \eqref{eq:genera} give
\[
 d\le\frac{(q^2-1)Q-2}{(q-1)(N-\delta)-2}
 <\frac{q^2-1}{q-1-q^{-1}}<q+2,
\]
where the last inequality follows from $q^2>2(q+1)$ for $q\ge3$.
Thus $d\le q+1$.

If $r\le n-3$, the function $z$ has degree $q^r$ on $T$ by
Lemma~\ref{lem:ramification}. Its pullback has degree $D=dq^r$
and $2qD\le2q(q+1)q^{n-3}<Q-1$, since $q^2>2(q+1)$ for $q\ge3$.
Lemma~\ref{lem:small} gives $K\mid dq^r$, hence $q^2-1\mid d$,
contrary to $1\le d\le q+1$.

It remains to consider $r=n-2,n-1$. Put $s=q^{n-r}$ and
\begin{equation}
h=\frac{z^{s(q-1)}}{1-y^{q-1}}.
\label{eq:h}
\end{equation}
Its only poles are the $\delta$ places over each
$y=\alpha\in\F_q^*$. By Lemma~\ref{lem:ramification} its order there
is $(s(q-1)R-N)/\delta=-(s+1)/\delta$.
Thus $\deg h=(q-1)(s+1)$ and
\[
D=\deg\pi^*h=d(q-1)(s+1)\le q^4-1.
\]
Since $2q(q^4-1)<q^6-1\le Q-1$, Lemma~\ref{lem:small} gives
$\pi^*h\in\overline{k}(w)$ and $q(q+1)\mid d(s+1)$.
Now $q\mid d$ because $\gcd(q,s+1)=1$, so $d=q$.
Of $s=q^2$ and $s=q$, only $s=q$ satisfies $q+1\mid s+1$.
Hence $r=n-1$ and $D=K$.

At this degree $\overline{k}(\pi^*h)=\overline{k}(w)$.
It follows that the target subfield $\pi^*\overline{k}(T)$ contains
$\overline{k}(w)$, so $F/\pi^*\overline{k}(T)$ is Galois.
Its group is a Sylow $p$-subgroup of $G_0$ and hence conjugate to
the translations. Conversely, Remark~\ref{rem:endpoint} identifies
the translation quotient over $k$ with $T_{q,n,n-1}$.
\end{proof}

\begin{proof}[Proof of Theorem~\ref{thm:main}\textup{(2)}]
For odd $n\ge5$ and $\delta=1$, maximality and \eqref{eq:genera} give
\begin{align*}
\#\X(k)-(q-1)\#T(k)&=qQ\{Q-q(q-1)\}-q+2>0,\\
(q+1)(2g(T)-2)-(2g(\X)-2)&=q(q-2)(q+1)>0 .
\end{align*}
Point counts and the Hurwitz genus formula therefore give $d=\deg\pi=q$.
Let $E=\pi^*\overline{k}(T)$ and $L=E\overline{k}(w)=F^H$, with
$H\le G_0$ of order $e\mid q$. After conjugation we may take
$H$ inside the translations. Each nonidentity translation
$\sigma_a$ fixes only the place $O$ above $u=\infty$. There
$\ord_O(u)=-N$, $\ord_O(Z)=-Q$, and $t=Z/u$ is a local parameter.
Since $\sigma_a(t)-t=-aZ/(u(u+a))$, its order is $Q+2$.
The different has degree $(e-1)(Q+2)$, so
\[
2e\,g(L)=(q^2-e)Q-q^2(q-1).
\]
If $e<q$, apply the Castelnuovo--Severi inequality to
$L=E\overline{k}(w)$
with indices $q/e,K/e$ and $2g(E)=(q-1)Q$.
Since $e\mid q$ and $e<q$, we have $e\le q/2$, so rearranging gives
\[
Q\le \frac{q^2(q-1)}{q-e}+\frac{2(K-e)}e<3q^3<q^5\le Q .
\]
Thus $e=q$ and $L=E$. The formula for $g(L)$ now gives
$2g(E)=(q-1)(Q-q)$, contradicting $2g(T)=(q-1)Q$.
\end{proof}

Changing the parameters of the covering curve gives further coverings.

If $d$ is a proper divisor of $n$, the curves $\X_{q^d,n/d}$ and $\X_{q,n}$ are maximal
over the same field $k$. Replacing $q,n$ by $q^d,n/d$ in the covering
curve gives the following quotients.

\begin{proposition}\label{prop:changed-source}
Let $d$ be a proper divisor of $n$ and put
$s=q^d$, $m=n/d$ and $\ell=(s-1)/(q-1)$. There is a Galois covering
over $k$
\begin{equation}
\X_{q^d,n/d}\longrightarrow T_{q,n,n-d}
\qquad\text{of degree }s\ell.
\label{eq:changed-source}
\end{equation}
\end{proposition}
\begin{proof}
Put $R_s=(q^{n-d}-1)/(s-1)$, so $R=\ell R_s$ for $r=n-d$.
Choose $\theta\in\F_{s^2}^*$ with $\theta^s=-\theta$ and use the
source model with $c=\eta^{-R_s}\theta$. This choice gives $c^Q=-c$. In odd characteristic $\theta^Q=(-1)^m\theta$ and
$R_s\equiv m-1\pmod2$, while in characteristic two $c\in\F_Q^*$.
The function $Y=\theta(u^s-u)$ satisfies
\[
 Z^N=\eta^{-R_s}Y(1-Y^{s-1})^{R_s}.
\]
Raising this equation to the power $\ell$ gives \eqref{eq:model} with $y=Y^\ell$ and $z=Z^\ell$, since $y^{q-1}=Y^{s-1}$.
As a polynomial in $u$, $y$ has degree $s\ell$.
Both $k(\X_{s,m})/k(u)$ and $k(y,z)/k(y)$ have degree $N$ by
the Kummer equations, so the covering has degree $s\ell$.
The $s\ell$ automorphisms
\[
 (u,Z)\longmapsto(\lambda^N u+b,Z\lambda),\qquad
 b\in\F_s,\quad\lambda^\ell=1,
\]
fix $y,z$, proving that the covering is Galois over $k$.
\end{proof}

For example, $d=2$ and $n=8$ give a covering of $T_{q,8,6}$
by $\X_{q^2,4}$ over $\F_{q^{16}}$, while Theorem~\ref{thm:main}(1)
excludes a covering by $\X_{q,8}$ over $\overline{k}$.

\begin{remark}\label{rem:coprime}
The targets in Proposition~\ref{prop:changed-source} have $r=n-d$
with $d\mid n$, so $\gcd(n,r)=d$. Thus for $\gcd(n,r)=1$ the only
covering supplied by that proposition is the degree-$q$ translation
quotient at $r=n-1$. In the subfamily with $r$ odd, this occurs
precisely when $n$ is even.

Suppose further that $q$ and $n\ge5$ are primes and $1<r<n$ is odd.
Then $T_{q,n,r}$ is not covered over $\overline{k}$ by any
$\X_{s,m}$ with $s^m=q^n$. Indeed, for a prime power $s$ and
$m\ge1$ this equality allows only $(s,m)=(q,n)$ and $(q^n,1)$.
In the first case, $r$ odd gives $\delta=1$ by Lemma~\ref{lem:arithmetic},
so Theorem~\ref{thm:main}\textup{(2)} applies. In the second case,
$\X_{q^n,1}$ is $k$-isomorphic to the Hermitian curve
$\mathcal H_{q^n}$ \cite{MW2026}, and $r\nmid n$ excludes a covering
by Theorem~\ref{thm:non-cover-by-Herm}.
\end{remark}

\section{Non-covering by the Hermitian curve}
\label{sec:noncover}

We prove the Hermitian covering criterion in
Theorem~\ref{thm:non-cover-by-Herm}.
By Proposition~\ref{prop:parameter-symmetry}, it suffices to consider
$1\le r\le n$. Proposition~\ref{prop:explicit-cover} settles the
case $r\mid n$. We prove that a covering
$\HH_Q\to T_{q,n,r}$ over $\overline{k}$ with $1<r<n$ forces
$r\mid n$.

\subsection{Polynomial representations and vanishing orders}
\label{sec:local}

The non-covering proof uses the two rational places
$P_0,P_\infty$ to restrict a hypothetical covering.
After reduction to a separable covering over $k$, their pullbacks
are $k$-rational divisors, and the genus bound makes the covering
degree less than $Q$. For functions with poles bounded by $N$
times either divisor, Lemma~\ref{lem:representation} gives ratios
of polynomials of degree less than $Q$;
Lemma~\ref{lem:local} then restricts their vanishing orders.
The single-pole functions constructed in Section~\ref{sec:poles}
then force the covering degree and the two fibres in
Propositions~\ref{prop:infinity} and~\ref{prop:zero}.
Combining these fibres with the Frobenius divisor relation gives
$y=h^R$ on $\HH_Q$; the defining equation then yields the
polynomial identity used to finish the proof in
Section~\ref{sec:completion}.

We use the model
\begin{equation}
\label{eq:Hermitian}
 \HH_Q:\quad b^Q+b=a^N,\qquad N=Q+1,
\end{equation}
whose unique place at infinity is denoted by $O$.
We compute functions and divisors over $\overline{k}$.
A divisor is $k$-rational if it is invariant under the
$Q^2$-power Frobenius.

For a nonzero homogeneous polynomial $F$ of degree $d<N$,
its intersection divisor on $\HH_Q$ is
\begin{equation}
\label{eq:intersection}
 E_F=(F(a,b,1))+dNO.
\end{equation}
The vanishing order of $F$ at $P$ is the coefficient of $P$ in
$E_F$. Since the equation of $\HH_Q$ is absolutely irreducible
of degree $N$, no nonzero polynomial of degree less than $N$
vanishes identically on $\HH_Q$.

\begin{lemma}
\label{lem:tangent}
Let $P$ be a geometric point of $\HH_Q$, and let $\Fr$ be the
$Q^2$-power Frobenius. The tangent line at $P$ has intersection
divisor $QP+\Fr(P)$. In particular, its intersection
multiplicity at $P$ is
\begin{equation}
\label{eq:kappa}
 \kappa(P)=
 \begin{cases}
 N,&P\in\HH_Q(k),\\
 Q,&P\notin\HH_Q(k).
 \end{cases}
\end{equation}
If $E$ is an effective $k$-rational divisor of degree $d$, then
$NE\sim dNO$.
\end{lemma}
\begin{proof}
For a finite point $P=(\alpha,\beta)$, set
$X=a-\alpha$ and $Y=b-\beta-\alpha^QX$.
Equation~\eqref{eq:Hermitian} becomes
\begin{equation}
\label{eq:tangent-local}
 Y^Q+Y=(\alpha-\alpha^{Q^2})X^Q+X^{Q+1},
\end{equation}
and the tangent line is $Y=0$.
Its intersection multiplicity at $P$ is $Q$, unless
$\alpha^{Q^2}=\alpha$, in which case it is $N$.
If $\alpha\in k$, the equality
$\beta^Q+\beta=\alpha^{Q+1}$ also gives $\beta^{Q^2}=\beta$.
In the other case, the residual intersection point is
$(\alpha^{Q^2},\beta^{Q^2})$.
The tangent line at $O$ meets the curve only at $O$, with
multiplicity $N$.

All line sections are linearly equivalent, so
$QP+\Fr(P)\sim NO$.
Summing this relation over $E$ and using $\Fr(E)=E$ gives
$NE\sim dNO$.
This is the Hermitian case of the Frobenius divisor relation
in \cite[Corollary~1.2]{FGT1997}.
\end{proof}

The following is the form of \cite[Lemma~3.1]{MW2026} used here.

\begin{lemma}
\label{lem:representation}
Let $E$ be an effective $k$-rational divisor on $\HH_Q$ of
degree $1\le d<Q$.
For nonzero functions $1,f_1,\ldots,f_m\in\LL(NE)$,
there are homogeneous polynomials $F_0,\ldots,F_m$ over
$\overline{k}$, all of degree $d$, such that
$f_i=F_i/F_0$ on $\HH_Q$ and
\begin{equation}
\label{eq:representation-divisors}
 E_{F_0}=NE,\qquad E_{F_i}=NE+(f_i).
\end{equation}
\end{lemma}
\begin{proof}
By Lemma~\ref{lem:tangent}, choose a function $h$ with
$(h)=NE-dNO$. Then $h,hf_i\in\LL(dNO)$.
The Weierstrass semigroup at $O$ is $\langle Q,N\rangle$:
the functions $a,b$ have pole orders $Q,N$, and this numerical
semigroup has $Q(Q-1)/2$ gaps, equal to the genus of $\HH_Q$.
A basis of $\LL(dNO)$ is
\[
 \{a^u b^v:0\le u<N,\ v\ge0,\ uQ+vN\le dN\}.
\]
Indeed, every element of $\langle Q,N\rangle$ can be written as
$uQ+vN$ with $0\le u<N$ and $v\ge0$, by reducing $u$ modulo $N$.
These representations are unique: equality of two such orders
gives $N\mid u-u'$, hence $u=u'$ and $v=v'$.
The monomials therefore have distinct pole orders and are linearly
independent. For any function in $\LL(dNO)$ with a pole at $O$,
subtracting a scalar multiple of the monomial with the same pole
order lowers that order. Repetition leaves a constant, proving
that the monomials $a^u b^v$ with $0\le u<N$, $v\ge0$ and
$uQ+vN\le dN$ span $\LL(dNO)$.

If $u+v>d$, then
$uQ+vN\ge(d+1)Q>dN$, since $d<Q$.
Thus every monomial in this basis has total degree at most $d$.
Homogenizing $h,hf_i$ to degree $d$ gives the required
polynomials. Equation~\eqref{eq:intersection} yields
$E_{F_0}=NE$ and $E_{F_i}=NE+(f_i)$.
\end{proof}

\begin{lemma}
\label{lem:local}
Let $F$ be a nonzero homogeneous polynomial of degree $d<Q$.
If its vanishing order at a geometric point $P$ of $\HH_Q$ is
$v$, then
\begin{equation}
\label{eq:local}
 v\bmod Q\le d.
\end{equation}
Here $v\bmod Q$ is the residue in $\{0,\ldots,Q-1\}$.
\end{lemma}
\begin{proof}
Let $L_P$ be the tangent line equation at $P$, and write
$F=L_P^mG$ with $L_P\nmid G$.
In affine coordinates centered at $P$ with tangent line $Y=0$,
$X$ is a local parameter and the curve has an expansion
$Y=\psi(X)$ of order $\kappa(P)\in\{Q,Q+1\}$.
The polynomial $G(X,0)$ is nonzero and has degree at most $d-m<Q$.
Every term containing $Y$ has order at least $Q$ after substitution,
so it cannot cancel the lowest-order term of $G(X,0)$.
Hence the vanishing order of $G$ is some $m'\le d-m$, and
$v=m\kappa(P)+m'$.
Its residue modulo $Q$ is $m'$ if $\kappa(P)=Q$, and $m+m'$
if $\kappa(P)=Q+1$; both are at most $d$.
\end{proof}

\begin{corollary}
\label{cor:cover-local}
Let $\pi:\HH_Q\to X$ be a nonconstant $k$-morphism of degree
$d<Q$, and let $P_*$ be a rational place of $X$.
Suppose that $f$ has its only pole at $P_*$, of order
$h\le N$ (we also allow $f=1$, $h=0$).
For every geometric place $P$ above $P_*$, with ramification
index $e=e(P\mid P_*)$, we have
\begin{equation}
\label{eq:cover-local}
 \bigl(e(N-h)\bmod Q\bigr)\le d.
\end{equation}
\end{corollary}
\begin{proof}
Apply Lemma~\ref{lem:representation} to $E=\pi^*P_*$.
The numerator representing $\pi^*f$ has degree $d$ and
vanishing order $e(N-h)$ at $P$.
Lemma~\ref{lem:local} gives \eqref{eq:cover-local}.
\end{proof}

\subsection{Functions with a prescribed pole}
\label{sec:poles}

Throughout this subsection and the next, assume $1<r<n$
and let $\pi:\HH_Q\to T_{q,n,r}$ be a separable $k$-morphism
of degree $\nu$. Put
\[
 A=q^r,\qquad s=q^{n-r},\qquad c=q-1,\qquad
 \ell=\left\lfloor\frac R2\right\rfloor+1.
\]
Thus $Q=sA$ and $N=csR+s+1$. Lemma~\ref{lem:arithmetic} gives
\begin{equation}
\label{eq:extension-bounds}
 s\ge q,\qquad R\ge q+1,\qquad
 \delta\le s+1,\qquad \delta^2\le Q,\qquad \delta<A.
\end{equation}
Indeed, if $\delta=q^h+1$ with $h=\gcd(n,r)$, then
$h\mid n-r$, $r/h\ge2$ and $n/h$ is odd and at least three.
Thus $(q^h+1)^2\le q^{3h}\le Q$; the other bounds follow
as well. The cases $\delta=1,2$ are immediate.

To determine the fibres above $P_0$ and $P_\infty$, we use
functions with poles only at one of these places and of order at
most $N$. The functions below are chosen so that the complementary
pole orders $N-h$ form arithmetic progressions.
Corollary~\ref{cor:cover-local} will then restrict the ramification
indices by bounding their residues modulo $Q$. Set
\begin{equation}
\label{eq:progression-functions}
 g_j=z^s\left(\frac{z^{cs}}{1-y^{q-1}}\right)^j,
 \qquad 0\le j\le\ell.
\end{equation}
Equation~\eqref{eq:divisors} gives
\begin{align}
\label{eq:power-poles}
 (z^i)_\infty&=iAP_\infty,&
 (z^i/y)_\infty&=(N-i)P_0 \qquad(0\le i\le s),\\
\label{eq:progression-poles}
 (g_j)_\infty&=(Q-cj)P_\infty,&
 (g_j/y)_\infty&=\bigl(N-s(1+cj)\bigr)P_0.
\end{align}
Indeed, at every place above $\alpha\in\F_q^*$ we have
$\ord(g_j)=(sR-j(s+1))/\delta$, whereas
$\ord_{P_0}(g_j)=s(1+cj)$ and
$\ord_{P_\infty}(g_j)=-Q+cj$.
The inequalities
\[
 sR-\ell(s+1)\ge0,\qquad \ell<R,\qquad q\ell>R,
 \qquad c\ell>R
\]
follow from $s\ge3$, $R\ge4$ and the choice of $\ell$.
For even $R$, the first is $(s-1)R\ge2(s+1)$;
for odd $R$ it is weaker.
Thus the pole divisors in \eqref{eq:power-poles} and
\eqref{eq:progression-poles} have the stated supports and orders.

The next lemma supplies the arithmetic steps in this argument:
the first turns bounded residues into an integer arithmetic
progression, and the second determines the ramification index
from the progression associated with $g_j$.

\Needspace{9\baselineskip}
\begin{lemma}
\label{lem:residues}
\begin{enumerate}
\item\label{lem:residues-progression}
Let $\kappa$ be a positive integer.
If integers $b_0,\ldots,b_h\in[0,D]$ have consecutive differences
congruent modulo $\kappa$ and $2D<\kappa$, then they form
an arithmetic progression in $\mathbb Z$.
\item\label{lem:residues-exclusion}
If an integer $\tau$ satisfies $1\le \tau\le R+1$ and
\[
 \tau(1+cj)\bmod A\le R+1,\qquad 0\le j\le\ell,
\]
then $\tau=R$.
\end{enumerate}
\end{lemma}
\begin{proof}
For the first assertion, the consecutive differences belong to
$[-D,D]$, whose length is less than $\kappa$.

For the second, $\tau=R+1$ gives residue $R+c>R+1$ at $j=1$,
so $\tau\le R$.
Consecutive differences of the residues are $c\tau$ or $c\tau-A$.
If both lay in $[-R-1,R+1]$, then
$(c-1)R\le c\tau\le R+1$, forcing $c=2$ and $R\le2\tau\le R+1$.
But then the residue at $j=1$ is $3\tau$, with $R+1<3\tau<A$.
Thus the differences are equal, and the residues form an integer
arithmetic progression starting at $\tau$.
If the common difference is positive, the final term is at least
$1+c\ell>R+1$. If it is negative and $\tau<R$, then
\[
 \ell(A-c\tau)\ge q\ell>R\ge \tau,
\]
so the final term would be negative. Hence $\tau=R$.
\end{proof}

\subsection{The two fibres}

\begin{proposition}
\label{prop:infinity}
The covering degree is $\nu=sR$, and
$\pi^*P_\infty=sR\,O_1$ for some $O_1\in\HH_Q(k)$.
\end{proposition}
\begin{proof}
For a geometric place $P\mid P_\infty$, write
$e=e(P\mid P_\infty)$.
Corollary~\ref{cor:cover-local}, applied to $1,z,\ldots,z^s$,
puts each residue $e(N-iA)\bmod Q$ in $[0,D]$.
Since $Q=sA$, these residues are all congruent to $e$ modulo $A$.
Write
\[
 \rho_i=e(N-iA)\bmod Q=(e\bmod A)+Am_i,\qquad 0\le i\le s.
\]
The integers $m_i$ lie in $[0,\lfloor D/A\rfloor]$
and have congruent consecutive differences modulo $s$.
By Lemma~\ref{lem:residues}\textup{(\ref{lem:residues-progression})}
they form an integer arithmetic progression.
Since $e\le\nu\le D<Q$, the endpoints are $\rho_0=\rho_s=e$.
Thus the progression is constant and $s\mid e$.

Write $e=s\tau$; then $\tau\le R+1$.
The residues given by $g_j$, divided by $s$, are
$\tau(1+cj)\bmod A\le D/s\le R+1$.
Lemma~\ref{lem:residues}\textup{(\ref{lem:residues-exclusion})}
gives $\tau=R$.
Thus every geometric place above $P_\infty$ has index $sR$.
Since their indices sum to $\nu<2sR$, there is exactly one.
It is fixed by Frobenius, and hence rational; its index is $\nu=sR$.
\end{proof}

\begin{proposition}
\label{prop:zero}
Every geometric place above $P_0$ has ramification index $R$.
In particular,
\[
 \pi^*P_0=R E_0,\qquad \deg E_0=s,
\]
where $E_0$ is effective and $k$-rational.
\end{proposition}
\begin{proof}
By Proposition~\ref{prop:infinity}, $\nu=sR$.
Since $c\ge2$, we have $2\nu=2sR\le csR=Q-s<Q$.
For $P\mid P_0$ with index $e$, Corollary~\ref{cor:cover-local}
applied to $z^i/y$ gives $ei\bmod Q\le\nu$ for $0\le i\le s$.
By Lemma~\ref{lem:residues}\textup{(\ref{lem:residues-progression})},
these residues form an integer arithmetic progression beginning
with $0,e$. Hence $se\le\nu$ and $e\le R$.
The residues given by $g_j/y$, divided by $s$, are
$e(1+cj)\bmod A\le R$.
Lemma~\ref{lem:residues}\textup{(\ref{lem:residues-exclusion})}
gives $e=R$.
The divisor $E_0=\pi^*P_0/R$ has degree $s$ and is $k$-rational
because $\pi$ and $P_0$ are defined over $k$.
\end{proof}

The two fibres are now determined: $P_\infty$ has a unique
rational place above it, with ramification index equal to the
covering degree $sR$, while $\pi^*P_0=RE_0$ with $\deg E_0=s$.
Together with $(y)=N(P_0-P_\infty)$, these identities give the
divisor of the pullback of $y$. The Frobenius divisor relation
will allow us to take its $R$-th root in the Hermitian function field.

\subsection{Completion of the proof}
\label{sec:completion}

\begin{proof}[Proof of Theorem~\ref{thm:non-cover-by-Herm}]
By Proposition~\ref{prop:parameter-symmetry}, we may assume
$1\le r\le n$. Then $\delta\le R\le(Q-1)/(q-1)$, so
\eqref{eq:genus} gives $2g(T_{q,n,r})\ge(q-2)Q+q>2$.
It remains to prove the necessity of $r\mid n$ for $1<r<n$.

Suppose that $1<r<n$ and that a covering over $\overline{k}$ exists.
We first reduce to a separable covering.
View $\overline{k}(T_{q,n,r})$ as a subfield of
$L_0=\overline{k}(a,b)$ through its pullback, and let $E$ be the
maximal intermediate field separable over $\overline{k}(T_{q,n,r})$.
Then $L_0/E$ is purely inseparable of degree $p^e$ for some $e\ge0$.
Since $\overline{k}$ is perfect,
\cite[Proposition~3.10.2(c)]{Stichtenoth2009} gives
$E=L_0^{p^e}=\overline{k}(a^{p^e},b^{p^e})$.
The Hermitian equation has coefficients in $\F_p$, so
$a\mapsto a^{p^e}$, $b\mapsto b^{p^e}$, fixing $\overline{k}$,
defines a $\overline{k}$-isomorphism $L_0\to E$.
Composing $\overline{k}(T_{q,n,r})\hookrightarrow E$ with its
inverse yields a separable covering $\pi:\HH_Q\to T_{q,n,r}$.
By Proposition~\ref{prop:descent}, this covering is defined over $k$.

Write $A=q^r$ and $s=q^{n-r}$, and identify functions on
$T_{q,n,r}$ over $\overline{k}$ with their pullbacks.
The group $\Aut(\HH_Q)$ acts transitively on $\HH_Q(k)$
\cite[Proposition~2.2]{GSX2000}.
After an automorphism of $\HH_Q$, the point $O_1$ above $P_\infty$
is the infinite place $O$ of \eqref{eq:Hermitian}.
Propositions~\ref{prop:infinity} and~\ref{prop:zero} give
\[
 (y)=RN(E_0-sO).
\]
As $E_0$ is $k$-rational, Lemma~\ref{lem:tangent} gives
$N(E_0-sO)\sim0$. Choose $h_0$ with $(h_0)=N(E_0-sO)$.
Then $(y/h_0^R)=0$, so $y/h_0^R$ is a nonzero constant.
Absorbing an $R$-th root of this constant into $h_0$ over
$\overline{k}$, we obtain
\begin{equation}
\label{eq:y-root}
 y=h^R,\qquad (h)_\infty=sNO.
\end{equation}
Substitution in \eqref{eq:model} gives
$z^N=\{\eta^{-1}(h-h^A)\}^R$.
Write $N=\delta m$ and $R=\delta S$, so that $\gcd(m,S)=1$.
After changing a constant, the last relation has the form
$z^m=f^S$.
For integers $u,v$ with $um+vS=1$, the function $f^u z^v$
has $S$-th power $z$ and $m$-th power $f$.
A further constant scaling gives a function $w$ with
\begin{equation}
\label{eq:common-root}
 w^m=h^A-h,\qquad z=\gamma w^S,\qquad
 (w)_\infty=\delta QO,\quad \gamma\in\overline{k}^*.
\end{equation}

Using the monomial basis in the proof of Lemma~\ref{lem:representation},
write $h=H(a,b)$ and $w=W(a,b)$, where $H,W$ are linear
combinations of monomials $a^u b^v$ with $0\le u<N$ and $v\ge0$.
Their pole orders $uQ+vN$ are distinct, so none can exceed the
pole order of the function it represents.
A monomial of total degree at least $s+1$ has pole order at least
$(s+1)Q>sN$; hence \eqref{eq:y-root} gives $\deg H\le s$.
Likewise, \eqref{eq:common-root} gives $\deg W\le\delta$.
Among these monomials, only $a^\delta$ has pole order $\delta Q$:
any monomial of degree $\delta$ involving $b$ has larger pole
order, whereas one of smaller degree has pole order at most
$(\delta-1)N<\delta Q$, since $\delta<Q$.
Thus $\deg W=\delta$, with highest homogeneous part
$\alpha a^\delta$ for some $\alpha\ne0$.

The function-field identity \eqref{eq:common-root} becomes
\begin{equation}
\label{eq:Hermitian-polynomial-identity}
 W^m=H^A-H+\lambda(a^N-b^Q-b),
 \qquad \lambda\in\overline{k}^*.
\end{equation}
The polynomial $W^m-H^A+H$ vanishes on the irreducible Hermitian
plane curve and has total degree $N$: its $a^N$ coefficient is nonzero,
whereas $\deg H^A\le sA=Q$. Thus the multiplier of the
Hermitian equation is a nonzero constant.

We show that $W$ is independent of $b$, so that comparison
of the $b$-coefficients in \eqref{eq:Hermitian-polynomial-identity}
will involve only $H$ and the Hermitian equation.
Indeed, if $d=\deg_b W>0$, then $d\le\delta<A$ and
$md\ge m=N/\delta>Q/A=s$.
In \eqref{eq:Hermitian-polynomial-identity}, every power of $b$
above degree $s$ on the right has exponent divisible by $A$,
since $A\mid Q$ and $\deg H\le s$.
The left side has $b$-degree $md$, so $A\mid md$.
As $\gcd(A,m)=1$, this forces $A\mid d$, a contradiction.

Write $H=\sum_i h_i(a)b^i$, taking $h_i=0$ for $i>\deg_b H$.
Since $W$ is independent of $b$ and
$H^A=\sum_i h_i(a)^A b^{Ai}$, comparison in
\eqref{eq:Hermitian-polynomial-identity} gives
\[
 h_1=-\lambda\ne0,\qquad
 h_{A^j}=h_{A^{j-1}}^A\ne0\quad(1\le j,\ A^j<Q).
\]
For the latter exponents, the only contributions on the right
are from $H^A$ and $-H$; the Hermitian term has $b$-degrees
only $0,1,Q$.
If $r\nmid n$, choose the first $j\ge1$ with $A^j>s=Q/A$.
Then $s<A^j<Q$, since equality at the upper endpoint would
give $rj=n$. This contradicts $\deg H\le s$.
Thus any covering requires $r\mid n$.

If $r\mid n$, including $r=1,n$, the required $k$-Galois covering
is given by Proposition~\ref{prop:explicit-cover}.
\end{proof}

\section{Automorphisms, isomorphisms and Galois subfields}
\label{sec:automorphisms}

Throughout this section assume \eqref{eq:range}, so $\delta=1$.
We determine the full automorphism group, classify the curves
$T_{q,n,r}$ up to isomorphism for fixed $q,n$, and give explicit
equations for all their Galois subfields.

Write $T=T_{q,n,r}$, $M=(q-1)N$, and $\mathcal A=\Aut(T)$.
The assignments
\begin{equation}
\label{eq:aut}
 \sigma_\lambda(y)=\lambda^Ny,\qquad
 \sigma_\lambda(z)=\lambda z,\qquad \lambda^M=1,
\end{equation}
preserve \eqref{eq:model}, since $\lambda^N\in\F_q^*$.
They form a cyclic subgroup $G\le \mathcal A$ of order $M$.
Since $M\mid Q^2-1$, all the multipliers $\lambda$ belong to $k$.
Let $C$ be the unique subgroup of $G$ of order $N$.
Then $\overline{k}(T)^C=\overline{k}(y)$.

\subsection{The full automorphism group}

\begin{theorem}
\label{thm:full-aut}
Assume \eqref{eq:range}. Then $\mathcal A=G$, so every
automorphism is one of the maps in \eqref{eq:aut}. In particular,
$\Aut(T_{q,n,r})\cong C_{(q-1)(q^n+1)}$.
\end{theorem}

We prove the theorem by determining first the normalizer of
$C$ and then the full stabilizers of its fixed points.
These stabilizers give three distinct short orbits of
$\mathcal A$ whose contributions to the Hurwitz genus formula
force $[\mathcal A:G]=1$.

Let $\mathcal B=\{P_\alpha:\alpha\in\F_q^*\}$.
By Lemma~\ref{lem:ramification},
\[
 \Fix(C)=\{P_0,P_\infty\}\cup\mathcal B.
\]
The three short $G$-orbits are
$\{P_0\}$, $\{P_\infty\}$ and $\mathcal B$, with stabilizers
$G,G,C$, respectively.
At every other geometric point, $z$ is finite and nonzero,
so \eqref{eq:aut} gives a trivial stabilizer.

\Needspace{5\baselineskip}
\begin{lemma}
\label{lem:normalizer}
The normalizers in $\mathcal A$ satisfy $\Nrm_{\mathcal A}(C)=G$ and
$\Nrm_{\mathcal A}(G)=G$.
\end{lemma}
\begin{proof}
Put $f(y)=\eta^{-R}y(1-y^{q-1})^R$.
If $\sigma\in\Nrm_{\mathcal A}(C)$, then its action on
$\overline{k}(y)$ is a fractional linear transformation
$\mu$.
The $C$-eigenspaces over $\overline{k}(y)$ are
$\overline{k}(y)z^i$, $0\le i<N$. Since $\sigma$ normalizes
$C$ and $z$ has a faithful character, we obtain
\begin{equation}
\label{eq:Kummer-lift}
 \sigma(z)=h(y)z^j,\qquad
 f(\mu(y))=h(y)^Nf(y)^j,
 \qquad \gcd(j,N)=1,
\end{equation}
for some $h(y)\in\overline{k}(y)^*$ and $1\le j<N$.
At the places of $\overline{k}(y)$ ramified in
$\overline{k}(T)/\overline{k}(y)$, the valuations of $f$ modulo $N$ are
\begin{equation}
\label{eq:branch-weights}
 1\quad(y=0),\qquad R\quad(y\in\F_q^*),\qquad
 -q^r\quad(y=\infty).
\end{equation}
Taking valuations in \eqref{eq:Kummer-lift} gives
\[
 \ord_{\mu(P)}(f)\equiv j\ord_P(f)\pmod N.
\]
Thus multiplication by $j$ permutes the residue classes of
these valuations modulo $N$, with their multiplicities. The three values in \eqref{eq:branch-weights} are distinct
modulo $N$ by \eqref{eq:arithmetic}.
Only $R$ occurs at more than one ramified place, namely at the
$q-1\ge2$ places in $\F_q^*$.
Consequently \eqref{eq:Kummer-lift} forces $\mu$ to
preserve $\F_q^*$ and gives $jR\equiv R\pmod N$.
As $\gcd(R,N)=1$ and $1\le j<N$, we have $j=1$.
The remaining two residues are distinct, so $\mu$ fixes $0$
and $\infty$ individually.
Hence $\mu(y)=ay$ with $a\in\F_q^*$.

Equation~\eqref{eq:Kummer-lift} now gives $h(y)^N=a$.
Hence $h(y)=\lambda\in\overline{k}^*$ with $\lambda^N=a$,
so $\sigma=\sigma_\lambda\in G$ and $\Nrm_{\mathcal A}(C)=G$.
Since $C$ is characteristic in $G$, we also have
$\Nrm_{\mathcal A}(G)\le\Nrm_{\mathcal A}(C)=G$.
\end{proof}

The first ramification group $\mathcal A_P^{(1)}$ is the unique Sylow
$p$-subgroup of $\mathcal A_P$, and $\mathcal A_P/\mathcal A_P^{(1)}$ is cyclic of order
prime to $p$; see \cite[Theorem~11.49]{HKT2008}.
The curve $T$ is maximal and hence has $p$-rank zero.
By \cite[Proposition~3.12]{GSY2015}, every nontrivial $p$-subgroup
of $\mathcal A$ fixes exactly one geometric point of $T$ and acts freely
elsewhere. We use these facts to determine $\mathcal A_P$ for
$P\in\Fix(C)$.

\begin{lemma}
\label{lem:stabilizers}
The stabilizers in $\mathcal A$ are
\begin{equation}
\label{eq:stabilizers}
 \mathcal A_{P_0}=\mathcal A_{P_\infty}=G,\qquad
 \mathcal A_{P_\alpha}=C\quad(\alpha\in\F_q^*).
\end{equation}
\end{lemma}
\begin{proof}
Let $P\in\Fix(C)$ and put $S=\mathcal A_P^{(1)}$.
Since $C\le \mathcal A_P$ and $S\trianglelefteq \mathcal A_P$, the group $C$
normalizes $S$.
Lemma~\ref{lem:normalizer} gives
\[
 \Cen_S(C)\le S\cap\Nrm_{\mathcal A}(C)=S\cap G=1.
\]
Suppose that $S\ne1$ and put $Y=T/S$. The group $C$ acts on
$Y$, and $S$ acts freely on $T\setminus\{P\}$.
The points of $\Fix(C)$ have distinct images on $Y$.
Indeed, the $S$-orbit of $P$ is a singleton. If
$U,V\in\Fix(C)\setminus\{P\}$ have the same image, then
$V=\tau(U)$ for some $\tau\in S$.
For every $c\in C$,
\[
 (c\tau c^{-1})(U)=V=\tau(U).
\]
Both $c\tau c^{-1}$ and $\tau$ belong to $S$. Since $S$ acts
freely at $U$, we get $c\tau c^{-1}=\tau$. Hence
$\tau\in\Cen_S(C)=1$, so $U=V$.

The induced action of $C$ on $Y$ is faithful.
Indeed, an element acting trivially on
$\overline{k}(T)^S$ belongs to
$\Gal(\overline{k}(T)/\overline{k}(T)^S)=S$, and
$C\cap S=1$.
Thus $C$ acts faithfully on $Y$ with at least $q+1$ fixed
points. The Hurwitz genus formula for $Y\to Y/C$ gives
\begin{equation}
\label{eq:wild-contradiction}
 2g(Y)-2\ge-2N+(q+1)(N-1)=2g(T)-2.
\end{equation}
But $T\to Y$ is separable of degree $|S|>1$, so
$2g(T)-2\ge|S|(2g(Y)-2)$.
This contradicts \eqref{eq:wild-contradiction}, since
$g(T)>1$.

Therefore $S=1$ and $\mathcal A_P$ is tame cyclic.
It contains $C$, hence normalizes $C$, and
Lemma~\ref{lem:normalizer} gives $\mathcal A_P\le G$.
The stabilizers $G_{P_0}=G_{P_\infty}=G$ and
$G_{P_\alpha}=C$ then give \eqref{eq:stabilizers}.
\end{proof}

\begin{proof}[Proof of Theorem~\ref{thm:full-aut}]
Put $h=[\mathcal A:G]$.
By Lemma~\ref{lem:stabilizers}, the $\mathcal A$-orbits of $P_0$ and
$P_\infty$ both have length $h$.
They are distinct: if $\sigma(P_0)=P_\infty$, then
$\sigma G\sigma^{-1}=G$, and hence
$\sigma\in\Nrm_{\mathcal A}(G)=G$ by Lemma~\ref{lem:normalizer}.
This is impossible because $G$ fixes these two points
individually.
The $\mathcal A$-orbit containing $\mathcal B$ has length $h(q-1)$;
its stabilizer has order $N<M$, so it is a third orbit.
The stabilizers on these three orbits are tame, so their
contributions to the different have degrees $h(M-1)$,
$h(M-1)$ and $h(q-1)(N-1)$.
Keeping these contributions in the Hurwitz genus formula for
$T\to T/\mathcal A$ and using $g(T/\mathcal A)\ge0$ gives
\begin{equation}
\label{eq:full-group-Hurwitz}
 2g(T)-2\ge-2hM+2h(M-1)+h(q-1)(N-1)
 =h(2g(T)-2).
\end{equation}
Since $g(T)>1$, we obtain $h=1$ and $\mathcal A=G$.
\end{proof}

\subsection{Isomorphism classes}

\begin{proof}[Proof of Theorem~\ref{thm:B}]
Theorem~\ref{thm:full-aut} gives the full automorphism group.
Fix $q,n$ and suppose that $T_{q,n,r}$ and $T_{q,n,r'}$
are isomorphic over $\overline{k}$.
An isomorphism conjugates their full automorphism groups and
therefore induces a fractional linear transformation of the
quotient coordinate $t$ in \eqref{eq:three-branch}.
The place $t=1$ has ramification index $N$, while
$0,\infty$ have index $M$.
Since $q>2$, the induced fractional linear transformation fixes
$1$ and either fixes $0,\infty$ individually or interchanges
them. It is therefore $t\mapsto t$ or $t\mapsto1/t$.

The exponent triples for the two cyclic extensions of degree
$M$ are
\begin{equation}
\label{eq:isomorphism-triples}
 (1,q^r-1,-q^r),\qquad
 (1,q^{r'}-1,-q^{r'}).
\end{equation}
An isomorphism multiplies all three exponents by one unit
modulo $M$, by the eigenspace argument used in
\eqref{eq:Kummer-lift}.
If the places $0,\infty$ are fixed individually, the exponent
$1$ forces that unit to be $1$, and hence
$q^r\equiv q^{r'}\pmod M$.
If they are interchanged, comparison at $0$ and $\infty$
instead gives $q^{r+r'}\equiv1\pmod M$.

The multiplicative order of $q$ modulo $N$ is $2n$:
it divides $2n$, whereas a proper divisor is at most $n$
and $0<q^e-1<N$ for $1\le e\le n$.
Reducing the preceding congruences modulo $N$ therefore gives
$r=r'$ in the first case and excludes the second, since
$0<r+r'<2n$.
The converse is immediate, proving \eqref{eq:intro-isomorphism}.
\end{proof}

\begin{proof}[Proof of Corollary~\ref{cor:C}]
For $n=2^m$, the quotient $r/\gcd(n,r)$ is odd.
Lemma~\ref{lem:arithmetic} gives $\delta=1$ for odd $r$,
and for every $1<r<n$ when $q$ is even.
There are $2^{m-1}-1$ odd parameters $1<r<n$, none dividing $n$.
For even $q$, removing the $m-1$ divisors $2,4,\ldots,2^{m-1}$
from the $2^m-2$ internal parameters leaves $2^m-m-1$.
Theorems~\ref{thm:A},~\ref{thm:non-cover-by-Herm} and~\ref{thm:B}
give the result.
\end{proof}

\subsection{Galois subfields}

Since $G=\Aut(T)$ is cyclic, it has a unique subgroup $G_d$
of order $d$ for each divisor $d\mid M$.
The function $t=y^{q-1}$ is fixed by $G$, and
$k(T)^G=k(t)$ by \eqref{eq:three-branch}.

\Needspace{11\baselineskip}
\begin{corollary}
\label{cor:quotients}
Let $d\mid M$ and put $L=M/d$.
The fixed field of $G_d$ is
\begin{equation}
\label{eq:fixed-field}
 k(T)^{G_d}=k(t,W),\qquad W=z^d,\qquad
 W^L=\eta^{1-q^r}t(1-t)^{q^r-1}.
\end{equation}
Its genus is
\begin{equation}
\label{eq:quotient-genus}
 g(T/G_d)=\frac{L-\gcd(L,q^r-1)}2
         =\frac{L-\gcd(L,q-1)}2.
\end{equation}
These fields are maximal over $k$. Every subfield $E$ containing
$k$ for which $k(T)/E$ is finite Galois occurs in this way.
\end{corollary}
\begin{proof}
The equation for $W=z^d$ in \eqref{eq:fixed-field} is
Eisenstein at $t=0$, so $[k(t,W):k(t)]=L$.
Since $t,W$ are fixed by $G_d$ and
$[k(T)^{G_d}:k(t)]=M/d=L$, this proves
\eqref{eq:fixed-field}.
The exponent triple at $0,1,\infty$ is
$1,q^r-1,-q^r$. Since $\gcd(L,q^r)=1$, the Hurwitz genus formula gives
\[
 2g(T/G_d)-2
 =-2L+2(L-1)+L-\gcd(L,q^r-1),
\]
also when $L=1$. This is the first expression in
\eqref{eq:quotient-genus}.
Since $\delta=1$ and $L\mid M$, equation~\eqref{eq:gcd-M}
gives $\gcd(L,q^r-1)=\gcd(L,q-1)$.
The quotients are maximal by Serre's theorem, and every
Galois subfield arises from a subgroup of the full group $G$.
\end{proof}

For fixed $q,n$, formula~\eqref{eq:quotient-genus} gives the same
genera of Galois subfields for all curves in Theorem~\ref{thm:B}.
Their fixed-field equations~\eqref{eq:fixed-field} retain the parameter $r$.

\raggedbottom
\Needspace{9\baselineskip}
\section*{Declaration on the use of AI}

The new families of curves were discovered through extensive searches
carried out using OpenAI's GPT-6 Astra.
Generative AI was also used to assist with mathematical calculations,
the development and checking of proofs, and the revision of the
manuscript. The text and reasoning were proofread and verified by the
authors. The authors take responsibility for the mathematical content
and the final text.

\Needspace{8\baselineskip}

\end{document}